\documentstyle[12pt]{article}
\renewcommand{\baselinestretch}{1.5}

\begin{document}
\begin{center}
\Large  Steinitz Class of Mordell Groups of  Elliptic Curves 
               \\ With Complex Multiplication \\ [1cm]

\normalsize Tong  LIU  \ \ \ \  and \ \ \ \ Xianke  ZHANG \\
    {\small ({Tsinghua University, Department of Mathematics,
Beijing 100084, P. R. China}) }
\end{center}
\date{}

\renewcommand{\baselinestretch}{1.5}
\renewcommand{\theequation}{\arabic{section}.\arabic{equation}}
\newcommand{\namelistlabel}[1]{\hfill\mbox{#1}}
\newenvironment{namelist}[1]{ %
\begin{list}{}
    {
      \let\makelabel\namelistlabel
      \settowidth{\labelwidth}{#1}
      \setlength{\leftmargin}{1.1\labelwidth}
    }
    }{ %
\end{list}}
\vskip 0.4cm
{\large Abstract.  } 
{Let $E$ be an elliptic curve having Complex Multiplication by the
full ring ${\cal O}_K$ of integers of $K={\bf Q}(\sqrt{-D})$, let
$H=K(j(E))$ be the Hilbert class field of
$K$. Then the Mordell-Weil group $E(H)$ is an ${\cal O}_K$-module,
and its Steinitz class St($E$) is studied.
 When $D$ is a prime number,
it is proved that $St(E)=1$ if
 $D\equiv 3\pmod{4}$;  and $St(E)=[{\cal P}]^t$
  if $p\equiv 1\pmod{4}$, where
$[{\cal P}]$ is the ideal class  of $K$ represented by  prime
 factor ${\cal P}$ of
2 in $K$, $t$ is a fixed integer. General structures are
also discussed for $St(E)$  and for  modules over Dedekind domain.
These results develop the results
by D. Dummit and W. Miller for $D=10$ and some elliptic curves to
more general $D$ and elliptic curves.}

\vskip 0.3cm


\section{introduction}\label{pg: 3.1}

\qquad Let $K={\bf Q}(\sqrt {-D})$ be an imaginary quadratic 
number field,
${\cal O}_K$ the ring of all integers of $K$.
Let $E$ be an elliptic curve having Complex Multiplication
 by the ring ${\cal O}_K$. $E$ is defined over the field  
$F={\bf Q}(j(E))$,
$j(E)$ is the $j$-invariant of $E$.
So $H=K(j(E))$ is the Hilbert class of $K ^{[3]}$,
and the Mordell-Weil group $E(H)$ of $H$-rational points of $E$
is then naturally a module over the Dedekind domain  ${\cal O}_K$ 
(via complex multiplication). 
By the structure theorem for finitely generated modules over
Dedekind domain we have that
$$ E(H)\cong E(H)_{tor}\oplus {\cal O}_K \oplus\cdots\oplus 
{\cal O}_K\oplus {\cal A}
=E(H)_{tor}\oplus {\cal O}_K^{s-1} \oplus {\cal A},$$
here ${\cal A}$ is an ideal of ${\cal O}_K$ which is uniquely
determined up to a multiplication  by a number from $K$.
Thus $E(H)$ uniquely defines an ideal class $[{\cal A}] $ 
of $K$ represented by ${\cal A}$ ,
which is said to be the Steinitz class of $E$ and denoted by $St(E)$.
(Similarly, any module $M$ over a Dedekind domain $R$ defines
an ideal class of $R$, which is said to be the Steinitz class of $M$
and denoted by $St(M)$.)
So the structure of the Mordell group $E(H)$, as a module over the 
Dedekind domain ${\cal O}_K$,is uniquely determined by its 
rank $s$, Steinitz class $St(E)$, and torsion part. Therefore, it
is important to determine the Steinitz class $St(E)$.
 D. Dummit and W. Miller $^{[5]}$ in 1996 determined the Steinitz class 
for some specific elliptic curves when $D=10$ and also found some properties  
of them. \par

   Let $l=rank_{\bf Z}(E(F))$ be the ${\bf Z}$-rank of $(E(F))$,
$G=Gal(H/F)$ be the Galois group of $H/F$ (which will be shown
to be a quadratic extension). 
Let $[{\cal A}]$ denote the ideal
class represented by the ideal ${\cal A}$ of ${\cal O}_K$.
Since $St(E)$ is concerned only with the free part of $E(H)$,
we put $E(\cdot)_f=E(\cdot)/E(\cdot)_{tor}$, i.e.,
 the quotient group of the Mordell
group $E(\cdot)$ modulo its torsion part. Note that $E(\cdot)_f$ 
is isomorphic to the free part of $E(\cdot)$. We will also use
this notation to subgroups of Mordell groups.

We will analysis the interior structure of $E(H)$,
give a general theorem for the structure of modules over 
Dedekind domain, and then determine Steinitz classes $St(E)$
for some types of elliptic curves. In particular,
   when $D=p$ is a prime number and $p\equiv 3\pmod{4}$,
we will prove that $St(E)$ is the principal class of $K$£»
And when the prime number $D=p\equiv 1\pmod{4}$, we will  show that
$$St(E)={[{\cal P}]}^t,\ \ \ \ \ t=l+\log{|H^1(G,E(H)_f)|}$$
where $|H^1(G,E(H)_f)|$ is the order of the first 
cohomology group $H^1(G,E(H)_f)$, and $\cal P$ is any prime factor
of 2 in $K$. 

The Weierstrass equation of  $E$ could be assumed as  
$^{[3]}$
$$E:  \quad y^2=f(x)=x^3+a_2 x^2+a_4 x+a_6$$
with $a_2,\  a_4, \ a_6 \in F $. \par

\setcounter{equation}{0}
\section{Structure of $E(H)$}\label{pg: 3.2}

\quad {\bf Lemma 1}. The degree of the extension $H/F$ is $[H: F]=2$.\par

{\bf Proof}.  Obviously  $[H: F] \leq 2$. Now if $[H: F]=1$,
 then $K \subset F$. Consider $F={\bf Q}(j(E))$,
where $E$ could be any elliptic curve having complex multiplication 
by ${\cal O}_K$. Since the complex multiplication domain
${\cal O}_K$ is a ${\bf Z}$-module of rank 2,
 so the $j-$invariant $j(E)$
is a real number $^{[6]}$.
Thus $F={\bf Q} (j(E))$ has a real embedding into the complex field.
Note that $K$ is totally imaginary, $K \subset F$ is impossible,
so $H\not\subset F$,
$[H: F]=[K(j(E)): {\bf Q}(j(E))]=2$. This proves the lemma.
\par \vskip 0.33cm

For any $ \alpha \in {\cal O}_K$, let $[\alpha]$ denote   
the endomorphism of $E$ corresponding to $\alpha $.
Comparing to $E$, we consider the following elliptic curve
$$E_D:  \hskip 0.4cm -Dy^2=f(x).$$
Note that $E_D$ and $E$ are isomorphic via the map
$$ i :  \hskip 0.6cm  E_D({\bf C})\rightarrow E({\bf C}),\qquad  
(x,y)\rightarrow (x,\sqrt {-D} y).$$
Thus $E_D$ also has complex multiplication by ${\cal O}_K$,
and is defined over $F$. Via the isomorphism $i$ of $E$ and $E_D$,
we have obviously that
$$E_D(F) \cong I=\{ (x,\sqrt {-D}y)|(x,\sqrt {-D}y)
\in E(H),\ x,y\in F\}\subset E(H).$$
The subgroup $I$ of
$E(H)$ defined here is very important in the following analysis.
\par

{\bf Lemma 2}.  The map $i\circ [\sqrt{-D}]$ is an $F$-isogeny of $E$ 
to $E_D$. Thus 
        $$rank_{\bf Z}(E_D(F))=rank_{\bf Z}(E(F))=l$$

{\bf Proof}.   By [5] we have 
$$[\sqrt{-D}](x,y)=(a(x),y\sqrt{-D}b(x)),$$
with $a(x),\ b(x)\in F(x)$.
So  $\ i\circ [\sqrt{-D}]\ $ is an $F$-isogeny of  $E$ to $E_D$.\par
 \vskip 0.36cm 

{\bf Lemma 3}.  $\ (I_f: [\sqrt{-D}]E(F)_f)(E(F)_f: [\sqrt{-D}]I_f)=D^l$
\par \vskip 0.3cm

{\bf Proof}.  
\begin{eqnarray*}
D^l &=&(E(F)_f: [D]E(F)_f) \\
&=&(E(F)_f: [\sqrt{-D}]I_f)([\sqrt{-D}]I_f: [D]E(F)_f)\\
&=&(E(F)_f: [\sqrt{-D}]I_f)(I_f:  [\sqrt{-D}]E(F)_f).
\end{eqnarray*}
\par

{\bf Lemma 4}.  \hskip 0.8cm 
$\ 2E(H)_f \subset E(F)_f\oplus I_f \subset E(H)_f$,
$$rank_{\bf Z}(E(H))=rank_{\bf Z}(E(F))+rank_{\bf Z}(E_D(F))=
2\ rank_{\bf Z}(E(F))=2l.$$
 \par \vskip 0.3cm
{\bf Proof}. \ If $ P=(x,y) \in E(F)_f$ with $P \in I_f$, then 
 $y=0$,  which
means that $P$ is a torsion point. So $P=O$ is the infinite point, and 
$E(F)_f\oplus I_f=E(F)_f+ I_f\subset E(H)_f$.
For any $ Q\in E(H)_f,$ we have $ 2Q=(Q+Q^{\sigma})+(Q-Q^{\sigma})$ ,
 where $\sigma
 \in G$.
Via the definition of $E(F)_f$ and $I_f$, we have 
$$E(F)_f=\{P|P^{\sigma}=P,\forall P\in 
E(H)_f\},\qquad  I_f=\{P|P^{\sigma}=-P,\forall P\in E(H)_f\}.$$
So  $\ Q+Q^{\sigma}\in E(F)_f,\: Q-Q^{\sigma}\in I_f\: $,
$2Q\in E(F)_f\oplus I_f$.
Thus  $2E(H)_f \subset E(F)_f\oplus I_f \subset E(H)_f$.
This completes the proof.\par \vskip 0.3cm 

As for the index of $E(F)_f\oplus I_f$ in $E(H)_f$, we have
the following theorem. \par \vskip 0.4cm 

{\bf Theorem 1}.   $$(E(H)_f: E(F)_f\oplus I_f)=
\frac{2^l}{|H^1(G,E(H)_f)|},$$
where $|H^1(G,E(H)_f)|$ is the order of the cohomology group 
$H^1(G,E(H)_f)$.\\  \par

{\bf Proof}.  Consider $H^1(G,E(H)_f)=Z^1(G,E(H)_f)/B^1(G,E(H)_f)$.
Let $\ T=\{ P-P^{\sigma}|P\in E(H)_f ,\  \sigma \in G\}$.
we will prove that $Z^1(G,E(H)_f)\cong I_f$,        
$B^1(G,E(H)_f)\cong T$.

For any cocycle $\xi \in Z^1(G,E(H)_f)$ ,
let $\xi \stackrel{\phi}{\rightarrow}{\xi}_{\sigma} $.
By the definition of cocycle we have that 
$0={\xi}_e={\xi}_{{\sigma}^2}={({\xi}_{\sigma})}^{\sigma}+{\xi}_{\sigma}$
,  so   ${({\xi}_{\sigma})}^{\sigma}=-{{\xi}_{\sigma}}$,   
thus   ${\xi}_{\sigma}\in I_f$, 
and  $\phi$ is a map of $Z^1(G,E(H)_f)$ to $I_f$. Via the
map $\phi$ we could see  that 
$Z^1(G,E(H)_f)\cong I_f,\ B^1(G,E(H)_f)\cong T$.
Now consider the homomorphism 
$E(H)_f\stackrel{\psi= P-P^{\sigma}}{\longrightarrow}{T}$. 
 Obviously $2I_f\subset T$. Since ${\psi}^{-1}(2I_f)=E(F)_f\oplus I_f$,
so $$(E(H)_f: E(F)_f\oplus I_f)=(T: 2I_f)=
{(I_f: 2I_f)}/{(I_f: T)},$$
$$={2^l}/{|H^1(G,E(H)_f)|}.$$

\setcounter{equation}{0}
\section{ Main Results and Their Proofs}\label{pg: 3.3}

\qquad We first give a general theorem on torsion-free 
finitely-generated module over Dedekind domain, 
which establishes a relationship between
the Steinitz class and the index of the module in its corresponding
free module. This theorem is the key to our final 
results about Steinitz  class.
\par \vskip 0.4cm

{\bf Theorem 2}. Suppose that $L$ is a free ${\cal O}_K$-module,
and $M\subset L$ is a submodule,  $(L: M)<+\infty$.
Then there is an integral ${\cal O}_K$-ideal ${\cal A}$ such that 
$[{\cal A}]$ is the Steinitz class of $M\ $,      
and $N^K_{\bf Q}({\cal A})=(L: M)\ $,
where $ \ N^K_{\bf Q}(\cdot)\ $ is the norm map of ideals from $K$ 
 to the rationals ${\bf Q}$.\par \vskip 0.3cm

{\bf Proof}.  Let $L=\bigoplus\limits_{i=1}^{n}{{\cal O}_K}e_i\ $, so
 $\ \{e_1,\cdots ,e_n\}$ is an ${\cal O}_K$-basis for $L$.
We will inductively prove that there is ${\cal O}_K$-ideals 
${\cal B}_i\ (i=1,\cdots ,n)$ such that
$M\cong \bigoplus\limits_{i=1}^{n}{{\cal B}_i}\  $,  
 and $(L: M)=\prod\limits_{i=1}^{n}({\cal O}_K: {\cal B}_i)$.

When $n=1$, every thing is obvious. Assume then the statement is true for 
$n-1$ and consider the homomorphism of ${\cal O}_K$-modules: 
$\rho : \ L\rightarrow {\cal O}_K$,\hskip 0.3cm 
$\rho (\sum\limits_{i=1}^{n}{{r_i}{e_i}})=r_n$.
Then ${\cal B}=\rho (M)$ is an ideal of ${\cal O}_K$,
 and the sequence  
$$0\rightarrow N\rightarrow M
\stackrel{\rho }{\rightarrow}{\cal B}\rightarrow 0$$
is exact,  where $N=ker(\rho )\cap M$.
Since ${\cal B}$ is a projective ${\cal O}_K$-module, there exists 
${\cal O}_K$-module ${\cal C}\subset M$ such that 
 ${\cal C}\cong {\cal B}$,  $\  \rho ({\cal C})={\cal B}$,
$\  M=N\oplus{\cal C}\cong N\oplus {\cal B}$. Thus 
$$(L: M)=(L: N\oplus{\cal C})=(L: \bigoplus\limits_{i=1}^{n-1}
{{\cal O}_K}+{\cal C})(\bigoplus\limits_{i=1}^{n-1}
{{\cal O}_K}+{\cal C}: N\oplus{\cal C})$$
where $(L: \bigoplus\limits_{i=1}^{n-1}
{{\cal O}_K}+{\cal C})=(\rho ^{-1}({\cal O}_K)): \rho ^{-1}({\cal B}))
=({\cal O}_K: {\cal B})$.\par

Consider ${\cal C}\cap \bigoplus\limits_{i=1}^{n-1}{{\cal O}_K}=
{\cal C}\cap ker(\rho )$. When restricted on  ${\cal C}$,
 the map $\rho $ is injective, so we have 
$$\bigoplus\limits_{i=1}^{n-1}{{\cal O}_K}+{\cal C}=
\bigoplus\limits_{i=1}^{n-1}{{\cal O}_K}\oplus {\cal C},$$
\begin{eqnarray*}
(\bigoplus\limits_{i=1}^{n-1}
{{\cal O}_K}+{\cal C}: N\oplus{\cal C})&=&
(\bigoplus\limits_{i=1}^{n-1}
{{\cal O}_K}\oplus {\cal C}: N\oplus{\cal C}) \\
&=&(\bigoplus\limits_{i=1}^{n-1}{{\cal O}_K}: N).
\end{eqnarray*}

Note that $N\subset \bigoplus\limits_{i=1}^{n-1}{{\cal O}_K}$.
So via the hypothesis of our induction, we know that there are 
${\cal O}_K$-ideals 
${\cal B}_i \ (i=1,\cdots ,n-1) $ such that 
$N\cong \bigoplus\limits_{i=1}^{n-1}{{\cal B}_i}$, and 
$(\bigoplus\limits_{i=1}^{n-1}{{\cal O}_K}: N)=\prod\limits_{i=1}^{n-1}
{({\cal O}_K: {\cal B}_i)}.$ Thus we have
$\; M\cong \bigoplus\limits_{i=1}^{n}{{\cal B}_i}\;  $ and 
$(L: M)=\prod\limits_{i=1}^{n}
{({\cal O}_K: {\cal B}_i)}=\prod\limits_{i=1}^{n}{N_{\bf Q}^K({\cal B}_i)}=
N_{\bf Q}^K(\prod\limits_{i=1}^{n}{{\cal B}_i})$,
where ${\cal B}_n={\cal B}$.
Now the proof is completed by the  following lemma. \par \vskip 0.36cm

{\bf Lemma 5}.  Assume ${\cal A}_1$ and ${\cal A}_2$ are two non-zero ideals of
Dedekind domain $R$,
then we have isomorphism of $R$-modules : 
${\cal A}_1\oplus{\cal A}_2\cong R\oplus 
{\cal A}_1{\cal A}_2$.\par
{\bf Proof}.  See Lemma 13 in [7] p.168.

We now intend to prove our main results via our Theorem 2. First
we need to find what the $L$ and $M$ of Theorem 2 correspond to 
in $E(H)$ . \par \vskip 0.36cm

{\bf Lemma 6}.   $L={\cal O}_K\cdot E(F)_f$ is 
a free ${\cal O}_K$-module of rank $l$.\par \vskip 0.3cm 

{\bf Proof}. Assume $P_1,\cdots ,P_l\ $ form a ${\bf Z}$-basis 
of  $E(F)_f$. We will prove  
$$L={\cal O}_K\cdot E(F)_f=\bigoplus\limits_{i=1}^{l}{{{\cal O}_K}P_i}.$$
Now suppose that $\sum\limits_{i=1}^{l}{{[{\alpha}_i]}P_i}=0$
for some ${{\alpha}_i} \in {\cal O}_K\ (i=1,\cdots , l).$  
When $D\equiv 3 \pmod{4}$, 
we have $ \alpha _i=
s_i + t_i (1+{\sqrt{-D}}/2$
$(s_i,t_i\in{\bf Z},\ i=1,\cdots ,l)$,
then via$\sum\limits_{i=1}^{l}{{[{\alpha}_i]}P_i}=0$
we have $\sum\limits_{i=1}^{l}{{[2s_i+t_i]}P_i}=0$ and 
$\sum\limits_{i=1}^{l}{{[\sqrt{-D}t_i]}P_i}=0$.
Thus $t_i=0,\: s_i=0,\: {\alpha}_i=0 \: (\ i=1,\cdots ,l)$.
This proves the theorem when $D\equiv 3 \pmod{4}$.
The case $D\equiv 1 \pmod{4}$ goes in the same way.\par \vskip 0.4cm

The corresponding free modules $F$ for $M$ varies for different $D$.
First consider the case $D\equiv 3 \pmod{4}$.\par \vskip 0.3cm

{\bf Theorem 3}.  For $D\equiv 3 \pmod{4}$, we have $|H^1(G,E(H)_f)|=1$,   
and $E(H)_f={\cal O}_K\cdot E(F)_f+ I_f$.\par \vskip 0.3cm

{\bf Proof}. Let $P_1,\cdots ,P_l$ form a ${\bf Z}$-basis of $E(F)_f$,
and $Q_1,\cdots ,Q_l$ form a  ${\bf Z}$-basis of $I_f$.
Put $\alpha=({1+\sqrt{-D}})/{2}$.
We need only to prove that $E(H)_f/(E(F)_f\oplus I_f)
=C_1\oplus \cdots  \oplus C_l$,
where $C_i=(\ \overline{[\alpha]P_i}\ )$ is subgroup of order 2 generated by
$\overline{[\alpha]P_i}$ in the quotient group 
$E(H)_f/(E(F)_f\oplus I_f)$.
Obviously we have $\overline{[\alpha]P_i}\neq\overline{0}$;
otherwise there would be $ t_j,\ s_j\in {\bf Z}
\ (j=1,\cdots ,l)$ such that
$[\alpha]P_i=\sum\limits_{j=1}^{l}{[t_j]P_j}+
\sum\limits_{j=1}^{l}{[s_j]Q_j}$,
then  $[1+\sqrt{-D}]P_i=\sum\limits_{j=1}^{l}{[2t_j]P_j}+
\sum\limits_{j=1}^{l}{[2s_j]Q_j}$,
and $P_i=\sum\limits_{j=1}^{l}{[2t_j]P_j}$,
 giving a contradiction. \par
Furthermore, if  $\sum\limits_{i=1}^{l}{[u_i]\overline{[\alpha]P_i}}=0$
for some $ u_i\in {\bf Z}\ (i=1,\cdots ,l)$,  then there are 
$ t_i,\ s_i\in {\bf Z}
\ (i=1,\cdots ,l)$ such that
$\sum\limits_{i=1}^{l}{[u_i\alpha]P_i}=
\sum\limits_{i=1}^{l}{[t_i]P_i}+\sum\limits_{i=1}^{l}{[s_i]Q_i}$,
   so 
$$\sum\limits_{i=1}^{l}{[u_i]P_i}+\sum\limits_{i=1}^{l}{[u_i\sqrt{-D}]P_i}=
\sum\limits_{i=1}^{l}{[2t_i]P_i}+\sum\limits_{i=1}^{l}{[2s_i]Q_i}.$$
Thus $\sum\limits_{i=1}^{l}{[u_i]P_i}=\sum\limits_{i=1}^{l}{[2t_i]P_i}$,
which gives $u_i=2t_i\ (i=1,\cdots ,l)$. Hence 
$[u_i]\overline{[{\alpha}]P_i}=\overline{[t_i][2{\alpha}]P_i}=
\overline{[t_i(1+\sqrt{-D})]P_i}=\bar 0$
This completes the proof.
\par \vskip 0.3cm 

Now we could prove our main results via Theorem 2.
\par \vskip 0.4cm

{\bf Theorem 4}. Suppose that $D=p\equiv 3 \pmod{4}$ is a prime number,
and elliptic curve $E$ has complex multiplication by the 
full ring ${\cal O}_K$
 of integers of $K={\bf Q}(\sqrt {-D})$. Then the  Steinitz class of 
$E$ is the principal class, i.e. $St(E)=1$.\par \vskip 0.3cm

{\bf Proof}.  Let $L={\cal O}_K \cdot E(F)_f$, $M=[\sqrt{-p}]E(H)_f$.
Since $M\cong E(H)_f$, so we need only to prove $St(M)$ 
is the principal class.\par

By Theorem 3 we have $E(H)_f={\cal O}_K\cdot E(F)_f
+ I_f$. Thus $$M=[\sqrt{-p}]E(H)_f=E(F)_f\cdot (\sqrt{-p}{\cal O}_K)+
[\sqrt{-p}]I_f\subset {\cal O}_K \cdot E(F)_f=L; $$
\begin{eqnarray*}
(L: M)&=&({\cal O}_K\cdot E(F)_f: [\sqrt{-p}]E(H)_f)\\
&=&\frac{(E(H)_f: [\sqrt{-p}]E(H)_f)}{(E(H)_f: {\cal O}_K \cdot E(F)_f)}\\
&=&\frac{p^l}{(E(H)_f: {\cal O}_K \cdot E(F)_f)}.
\end{eqnarray*}
Since $p$ is a prime number, so there is $t\ (0\leq t\leq l)$ such that 
$(L: M)=p^t$.
By Theorem 2, the Steinitz class of $M$ is equal to  $[{\cal A}]$ 
for some ${\cal O}_K$-ideal ${\cal A}$,
and $p^t=(L: M)=N_{\bf Q}^K({\cal A})$.  Since $p$ is a prime number,
${\cal A}={(\sqrt{-p}{\cal O}_K)}^t$ is principal.  Thus $St(E)=St(M)$
is the principal class.
\par \vskip 0.4cm

{\bf Theorem 5}. Suppose that $D=p\equiv 3 \pmod{4}$ is a prime number,
and  $E$ is an elliptic curve having complex multiplication by 
the ring ${\cal O}_K$
 of all integers of $K={\bf Q}(\sqrt {-D})$. Then the Steinitz class of 
$E$ is   $\ \ St(E)=[{\cal P}^t]$, \hskip 0.3cm 
 where $[{\cal P}]$ is the ideal class of $K$  represented
by ${\cal P}$ the prime factor of 2 in ${\cal O}_K$, 
  $\  2^t=2^l|H^1(G,E(H)_f)|$.\hskip 0.3cm 
In particular, the parity of $t$ determines $St(E)$ since ${\cal P}$ is not
principal while ${\cal P} ^2=2{\cal O}_K$ is.
\par \vskip 0.3cm

{\bf Proof}.  Let $L={\cal O}_K \cdot E(F)_f$,  $M=[2\sqrt{-p}]E(H)_f$. Since 
 $M\cong E(H)_f$,  so $St(E)=St(M)$.  
Note that $[2\sqrt{-p}]E(H)_f\subset [\sqrt{-p}](E(F)_f\oplus I_f)$,   
$[\sqrt{-p}]I_f\subset E(F)_f$.  Thus we have 
$M\subset {\cal O}_K \cdot E(F)_f=L$,  and 
\begin{eqnarray*}
(L: M)&=&( {\cal O}_K \cdot E(F)_f: [2\sqrt{-p}]E(H)_f)\\
&=&\frac{(E(H)_f: [2\sqrt{-p}]E(H)_f)}{(E(H)_f: {\cal O}_K \cdot E(F)_f)}\\
&=&\frac{{(4p)}^l}{(E(H)_f: E(F)_f\oplus I_f)(E(F)_f\oplus I_f: 
{\cal O}_K \cdot E(F)_f)}\\
&=&\frac{{(4p)}^l}{2^l |H^1(G,E(H)_f)|^{-1} (I_f: [\sqrt{-p}]E(F)_f)}\\
&=&2^l|H^1(G,E(H)_f)|\cdot {p^l}/{(I_f: [\sqrt{-p}]E(F)_f)}.
\end{eqnarray*}
Thus $(L: M)=2^tp^r$ for some $ t,r\geq 0$ since $p$ is a prime number.
By Theorem 2 we know that $N_{\bf Q}^K({\cal A})=2^tp^r$ for some 
${\cal O}_K$-ideal ${\cal A}$. Therefore 
${\cal A}={\cal P}^t{([\sqrt{-p}]{\cal O}_K)}^r$, $St(E)=[{\cal A}]=
[{\cal P}^t]$.   This proves the theorem.\par \vskip 0.35cm

{\bf {Corollary 1}}.  Suppose as in  Theorem 5.  
If $l=rank_{\bf Z}(E(F))=1$
 then $H^1(G,E(H)_f)$ 
determines the Steinitz class of $E$.\par \vskip 0.3cm 

Now we analysis the examples of Dummit and Miller in [5] by utilizing
our above method .
For these examples, we have  
$K={\bf Q}(\sqrt{-10})$,   
$\  D=10$,   $\  H=K(\sqrt{5})={\bf Q}(\sqrt{-10},\ \sqrt{5})$.
 We consider the ${\cal O}_K$-module  $L={\cal O}_K\cdot E(F)_f\  $ and 
$\  M=2[\sqrt{-10}]E(H)_f$. Then 
\begin{eqnarray*}
(L: M)&=&\frac{(E(H)_f: 2[\sqrt{-10}]E(H)_f)}{(E(H)_f: {\cal O}_K\cdot E(F)_f)}\\
&=&\frac{{(4\cdot 10)}^l}{(E(H)_f: E(F)_f\oplus I_f)(E(F)_f\oplus I_f: 
{\cal O}_K \cdot E(F)_f)}\\
&=&\frac{{(40)}^l}{{2^l}{|H^1(G,E(H)_f)|^{-1}}(I_f: [\sqrt{-10}]E(F)_f)}\\
&=&2^l|H^1(G,E(H)_f)|{10^l}/{(I_f: [\sqrt{-10}]E(F)_f)}.
\end{eqnarray*}
Thus the Steinitz class of $E$ is determined by the 2-exponent of
$2^l|H^1(G,E(H)_f)|(I_f: [\sqrt{-10}]E(F)_f)$. 

(DM1) Consider the following elliptic curve of Dummit and Miller: 
$$E_1:\hskip 0.3cm y^2=x^3+(6+6\sqrt 5 )x^2+(7-3\sqrt 5 ).$$ 
Then $l=1\ $, $|H^1(G,E(H)_f)|=2\ $, $(I_f: [\sqrt{-10}]E(F)_f)=1\ $,  
so  $2^l|H^1(G,E(H)_f)|(I_f: [\sqrt{-10}]E(F)_f)=4$.
Thus the Steinitz class of $E_1$ is the principal class, i.e.,
$St(E_1)=1$.

(DM2) Consider the elliptic curve   
$E_{1,isog}:\hskip 0.3cm y^2=x^3-912+12\sqrt 5 )x^2+(188+84\sqrt 5 )x $ 
in [5]. We have
 $l=1$,  $\  |H^1(G,E(H)_f)|=2$,  $\ (I_f: [\sqrt{-10}]E(F)_f)=2$,
 $\ 2^l|H^1(G,E(H)_f)|(I_f: [\sqrt{-10}]E(F)_f)=2^3$.
Thus the Steinitz class  $St(E_{1,isog})=[ {\cal P}]$,
where ${\cal P}$ is a prime factor of 2 in ${\cal O}_K$.

(DM3) For $E_{3}:\hskip 0.3cm y^2=x^3+612x^2+(46818-20808 \sqrt 5 )x $
 in [5], 
we have  $\ l=2$, $\ |H^1(G,E(H)_f)|=2$,  $\ (I_f: [\sqrt{-10}]E(F)_f)=1,\:  
2^l|H^1(G,E(H)_f)|(I_f: [\sqrt{-10}]E(F)_f)=2^3.$
   Thus $St(E_{3})= [{\cal P}]$, $\  {\cal P}$  
a prime factor of 2 in ${\cal O}_K$.
\par  \vskip 0.38cm  

There are still many open problems about the Steinitz classes of elliptic curves
having complex multiplication. For example, in the case 
$K={\bf Q}(\sqrt {-D})$ with prime $D\equiv 1 (\pmod{4})$,
we have the following conjecture.
\par \vskip 0.3cm

{\bf Conjecture}.   Both the cases $St(E)=1$ and $St(E)\not= 1$ exist
for some elliptic curves $E$ having complex multiplication by ${\cal O}_K$
with $K={\bf Q}(\sqrt {-D})$  mentioned above.\\
\par 

\section*{\large References }

\begin{namelist}{9999}
{\footnotesize
\item[1] Keqin Feng,  Introduction to Commutative Algebra,
Higher Education Press, Beijing,1985.
\item[2] Silverman J.H. , The Arithmetic of Elliptic Curves,
Spring-Verlag, New York, 1982.
\item[3]Shimura G.  , Introduction to The Arithmetic Theory of Automorphic Functions,Pul. Math. Soc. Japan11\\
(1971).
\item[4] Neukirch J. , Class Field Theory, Spring-Verlag, Berlin, New York,1980.
\item[5] Dummit D.S.  Miller W.L.  , The Steinitz Class of the Mordell Weil group of Some CM
    Elliptic Curves, J. of Number Theory 56(1996): 52-78.
\item[6]Gross  B.,  Arithmetic on Elliptic Curves with Complex Multiplication,
 SLN, Vol. 776,
Springer-Verlag, Berlin, 1980.}

\end{namelist}
\vskip 0.4cm 

Email:  xianke@tsinghua.edu.cn

(This paper was published in : Pacific Journal of Mathematics,
     193(2000), No.2, 371-379.   (IDS No.:  302KF)

available also  via
        http://nyjm.albany.edu:8000/PacJ/2000/193-2-7.html
)

\end{document}